УДК 519.856

# Двойственные методы поиска равновесий в смешанных моделях распределения потоков в больших транспортных сетях

*Гасников Александр Владимирович*[1,2,5] gasnikov@yandex.ru
*Гасникова Евгения Владимировна*[3] egasnikova@yandex.ru
*Нестеров Юрий Евгеньевич*[4,5] yurii.nesterov@uclouvain.be

[1] Кафедра Математических основ управления
Национального исследовательского Университета «Московский физико-технический институт».
141700, Россия, Московская область, г. Долгопрудный, Институтский переулок, д. 9
[2] Институт проблем передачи информации им. А.А. Харкевича Российской академии наук.
127051, Россия, г. Москва, Большой Каретный переулок, д.19 стр. 1
[3] Лаборатория структурных методов анализа данных в предсказательном моделировании (ПреМоЛаб)
Национального исследовательского Университета «Московский физико-технический институт».
141700, Россия, Московская область, г. Долгопрудный, Институтский переулок, д. 9
[4] Center for Operation Research and Econometrics Université Catholique de Louvain.
Voie du Roman Pays 34, L1.03.01 - B-1348 Louvain-la-Neuve (Belgium)
[5] Департамент больших данных и информационного поиска, факультет Компьютерных наук
Национального исследовательского Университета «Высшая школа экономики».
125319, Россия, г. Москва, Кочновский проезд, 3

**Аннотация**

В работе изучается задача поиска равновесного распределения потоков в транспортной сети, часть ребер которой характеризуется функциями затрат, а часть ребер характеризуется пропускной способностью и постоянными затратами на прохождение в отсутствии затора. Такие модели (получившие название *смешанные модели*) возникают, например, при описании грузоперевозок РЖД. Частным случаем смешанной модели является семейство моделей равновесного распределения потоков по путям: BMW-модель (модель Бэкмана), модель стабильной динамики. Поиск равновесия в смешанной модели сводится к решению задачи выпуклой оптимизации. В данной статье строится двойственная задача, которая решается методом зеркального спуска (двойственных усреднений). Приводятся два различных способа восстановления решения исходной (прямой) задачи. Показывается, что предложенные подходы допускают эффективное распараллеливание. Результаты о скорости сходимости предложенных численных методов соответствуют известным нижним оракульным оценкам для данного класса задач (с точностью до мультипликативных констант).

## 1. Введение

Поиск равновесий в транспортных сетях часто удается свести к поиску равновесия в соответствующей популяционной игре загрузки. Поиск равновесия (Нэша) в такой игре, в свою очередь, всегда сводится к решению задач оптимизации. При естественных предположениях (имеющих место на практике) такие задачи оказываются выпуклыми. Пожалуй, самым популярным примером тут является модель равновесного распределения потоков по путям BMW-модель [1, 2] (модель Бэкмана). В последнее время большую популярность приобрели модель стабильной динамики [1, 3, 4] и смешанная модель, объединяющая в себе свойства двух предыдущих моделей [5].

В данной работе результат п. 3 статьи [1] распространяются на смешанные модели. При этом также предлагается другая запись подхода п. 3 статьи [1], которая делает его более удобным для применения на практике.

Структура статьи следующая. В п. 2 описывается постановка задачи. Формулируется задачи выпуклой оптимизации, решение которой дает равновесие. Строится двойственная задача. Приводятся формулы связи прямых и двойственных переменных. В п. 3 описывается композитный вариант известного метода зеркального спуска (двойственных усреднений) для решения двойственной задачи. Приводится первый способ численного восстановления решения прямой задачи из (приближенного) решения двойственной задачи, полученного методом зеркального спуска. Доказывается теорема о скорости сходимости предложенного подхода, исходя из заданной точности восстановления решения прямой и двойственной задачи одновременно. В п. 4 приводится второй способ численного восстановления решения прямой задачи из (приближенного) решения двойственной задачи, полученного методом зеркального спуска. Доказывается аналогичная п. 3 теорема о скорости сходимости предложенного подхода, исходя из заданной точности восстановления решения прямой и двойственной задачи одновременно. В п. 5 проводится анализ (обсуждение) результатов, полученных в предыдущих пунктах. В самом конце приводится краткое резюме результатов численных экспериментов.

## 2. Описание задачи поиска равновесного распределения потоков в смешанной модели

Пусть транспортная сеть города представлена ориентированным графом $\Gamma=(V,E)$, где $V$ – узлы сети (вершины), $E\subset V\times V$ – дуги сети (рёбра графа), $O\subseteq V$ – источники корреспонденций ($S=|O|$), $D\subseteq V$ – стоки. В современных моделях равновесного распределения потоков в крупном мегаполисе число узлов графа транспортной сети обычно выбирают порядка $n=|V|\sim 10^4$ [1]. Число ребер $|E|$ получается в три четыре раза больше. Пусть $W\subseteq\{w=(i,j): i\in O, j\in D\}$ – множество корреспонденций, т.е. возможных пар «источник» – «сток»; $p=\{v_1,v_2,...,v_m\}$ – путь из $v_1$ в $v_m$, если $(v_k,v_{k+1})\in E$, $k=1,...,m-1$, $m>1$; $P_w$ – множество путей, отвечающих корреспонденции $w\in W$, то есть если $w=(i,j)$, то $P_w$ – множество путей, начинающихся в вершине $i$ и заканчивающихся в $j$; $P=\bigcup_{w\in W}P_w$ – совокупность всех путей в сети $\Gamma$; $x_p$ [автомобилей/час] – величина потока по пути $p$, $x=\{x_p:\ p\in P\}$; $f_e$ [автомобилей/час] – величина потока по дуге $e$:

$$f_e(x)=\sum_{p\in P}\delta_{ep}x_p \quad (f=\Theta x), \text{ где } \delta_{ep}=\begin{cases}1, & e\in p\\ 0, & e\notin p\end{cases};$$

$\tau_e(f_e)$ – удельные затраты на проезд по дуге $e$. Как правило, предполагают, что это – (строго) возрастающие, гладкие функции от $f_e$. Точнее говоря, под $\tau_e(f_e)$ правильнее понимать представление пользователей транспортной сети об оценке собственных затрат (обычно временных в случае личного транспорта и комфортности пути (с учетом времени в пути) в случае общественного транспорта) при прохождении дуги $e$, если поток желающих проехать по этой дуге будет $f_e$.

Рассмотрим теперь $G_p(x)$ – затраты временные или финансовые на проезд по пути $p$. Естественно считать, что $G_p(x) = \sum_{e\in E} \tau_e(f_e(x))\delta_{ep}$.

Пусть также известно, сколько перемещений в единицу времени $d_w$ осуществляется согласно корреспонденции $w \in W$. Тогда вектор $x$, характеризующий распределение потоков, должен лежать в допустимом множестве:

$$X = \left\{ x \ge 0: \sum_{p\in P_w} x_p = d_w,\, w \in W \right\}.$$

Рассмотрим игру, в которой каждой корреспонденции $w \in W$ соответствует свой, достаточно большой, набор однотипных “игроков”, осуществляющих передвижение согласно корреспонденции $w$ (относительный масштаб характеризуется числами $d_w$) Чистыми стратегиями игрока служат пути, а выигрышем – величина $-G_p(x)$. Игрок “выбирает” путь следования $p \in P_w$, при этом, делая выбор, он пренебрегает тем, что от его выбора также “немного” зависят $|P_w|$ компонент вектора $x$ и, следовательно, сам выигрыш $-G_p(x)$. Можно показать (см., например, [2]), что поиск равновесия Нэша–Вардропа $x^* \in X$ (макро описание равновесия) равносилен решению задачи

$$\Psi(f) = \sum_{e\in E} \sigma_e(f_e) = \sum_{e\in E} \int_0^{f_e} \tau_e(z)\,dz \to \min_{\substack{f=\Theta x\\ x\in X}}. \tag{1}$$

В пределе модели стабильной динамики (Нестеров–деПальма) [3, 4] на части ребер $E' \subseteq E$ может быть сделан следующий предельный переход

$$\tau_e^{\mu}(f_e) \xrightarrow[\mu\to 0+]{} \begin{cases} \bar{t}_e, & 0 \le f_e < \bar{f}_e \\ [\bar{t}_e, \infty), & f_e = \bar{f}_e \end{cases},$$

$$d\tau_e^{\mu}(f_e)/df_e \xrightarrow[\mu\to 0+]{} 0,\ 0 \le f_e < \bar{f}_e.$$

Задача (1) примет вид

$$\Psi(f) = \sum_{e\in E} \sigma_e(f_e) \to \min_{\substack{f=\Theta x,\, x\in X\\ f_e \le \bar{f}_e,\, e\in E'}},\ \sigma_e(f_e) = f_e \bar{t}_e,\ \ e \in E'.$$

В данной работе предлагается хорошо параллелизуемый двойственный численный метод поиска равновесия в смешанной модели (1), т.е. в модели, в которой часть ребер являются Бэкмановскими, а часть – Нестерова–деПальмы. Такие задачи возникают, например, в однопродуктовой модели грузоперевозок РЖД [5]. К таким “смешанным” моделям классический метод условного градиента (Франк–Вульфа), который зашит практически во все современные пакеты транспортного моделирования, к сожалению, уже не применим. Требуется разработка новых подходов.

Для задачи (1) можно построить следующую двойственную задачу [1, 4]

$$\Upsilon(t) = \underbrace{-\sum_{w\in W} d_w T_w(t)}_{F(t)} + \sum_{e\in E} \sigma_e^*(t_e) \to \min_{\substack{t_e \ge \bar{t}_e,\, e\in E'\\ t_e \in \operatorname{dom}\sigma_e^*(t_e),\, e\in E\setminus E'}}, \tag{2}$$

где $T_w(t)=\min_{p\in P_w}\sum_{e\in E}\delta_{ep}t_e$ – длина кратчайшего пути из $i$ в $j$ ($w=(i,j)\in W$) на графе $\Gamma$, ребра которого взвешены вектором $t=\{t_e\}_{e\in E}$. При этом решение задачи (1) $f$ можно получить из формул: $f_e=\bar{f}_e-s_e$, $e\in E'$, где $s_e\ge 0$ – множитель Лагранжа к ограничению $t_e\ge\bar{t}_e$; $\tau_e(f_e)=t_e$, $e\in E\setminus E'$. Заметим, что для ребер $e\in E'$ имеем $\sigma_e^*(t_e)=\bar{f}_e\cdot(t_e-\bar{t}_e)$, а для BPR-функций

$$\tau_e(f_e)=\bar{t}_e\cdot\left(1+\gamma\cdot\left(\frac{f_e}{\bar{f}_e}\right)^{\frac{1}{\mu}}\right)\Rightarrow\sigma_e^*(t_e)=\bar{f}_e\cdot\left(\frac{t_e-\bar{t}_e}{\bar{t}_e\cdot\gamma}\right)^{\mu}\frac{(t_e-\bar{t}_e)}{1+\mu}.$$

В приложениях обычно выбирают $\mu=1/4$. В этом случае приводимый ниже шаг итерационного метода (3) может быть осуществлен по явным формулам, поскольку существуют квадратурные формулы (формулы Кардано–Декарта–Эйлера–Феррари [6]) для уравнений 4-й степени.

Поиск вектора $t$ представляет самостоятельный интерес, поскольку этот вектор описывает затраты на ребрах графа транспортной сети. Решение задачи (2) дает вектор затрат $t$ в равновесии.

## 3. Первый способ восстановления решения прямой задачи по решению двойственной

Для решения двойственной задачи (2) воспользуемся методом зеркального спуска в композитном варианте [7–9] ($k=0,\dots,N$, $t^0=\bar{t}$, ограничение $t_e\in\operatorname{dom}\sigma_e^*(t_e), e\in E\setminus E'$ всегда будет не активным, т.е. его можно не учитывать)

$$t^{k+1}=\arg\min_{\substack{t_e\ge\bar{t}_e,\, e\in E'\\ t_e\in\operatorname{dom}\sigma_e^*(t_e),\, e\in E\setminus E'}}\left\{\gamma_k\left\{\left\langle\partial F(t^k),t-t^k\right\rangle+\sum_{e\in E}\sigma_e^*(t_e)\right\}+\frac{1}{2}\left\|t-t^k\right\|_2^2\right\}, \tag{3}$$

где $\partial F(t^k)$ – произвольный элемент субдифференциала выпуклой функции $F(t^k)$ в точке $t^k$, а

$$\gamma_k=\varepsilon/M_k^2,\ M_k=\left\|\partial F(t^k)\right\|_2,$$

где $\varepsilon>0$ – желаемая точность решения задач (1) и (2), см. (6). К сожалению, доказать неравенство (*) (см. доказательство теоремы 1 ниже) при $\mu>0$ получилось только при стратегии выбора шагов постоянными:

$$\gamma_k:=\varepsilon\Big/\left(\frac{1}{N+1}\sum_{k=0}^{N}M_k^2\right),\ \text{при этом}\ \tilde{M}_N^2:=\frac{1}{N+1}\sum_{k=0}^{N}M_k^2.$$

Положим

$$\bar{t}^N=\frac{1}{S_N}\sum_{k=0}^{N}\gamma_k t^k,\ S_N=\sum_{k=0}^{N}\gamma_k,$$

$$f_e^k\in-\partial_e F(t^k),\ \bar{f}_e^N=\frac{1}{S_N}\sum_{k=0}^{N}\gamma_k f_e^k,\ e\in E\setminus E', \tag{4}$$

$$\bar{f}_e^N=\bar{f}_e-s_e^N,\ e\in E', \tag{4'}$$

где $s_e^N$ – есть множитель Лагранжа к ограничению $t_e\ge\bar{t}_e$ в задаче

$$\frac{1}{S_N}\left\{\sum_{k=0}^{N}\gamma_k\cdot\left\{\sum_{e\in E'}\partial_e F\left(t^k\right)\cdot\left(t_e-t_e^k\right)\right\}+S_N\sum_{e\in E'}\bar{f}_e\cdot\left(t_e-\bar{t}_e\right)+\frac{1}{2}\sum_{e\in E'}\left(t_e-\bar{t}_e\right)^2\right\}\to\min_{t_e\ge\bar{t}_e,\,e\in E'}.$$

Критерий останова метода (правое неравенство)

$$\left(0\le\right)\Upsilon\left(\bar{t}^N\right)+\Psi\left(\bar{f}^N\right)\le\varepsilon. \tag{5}$$

В данной работе (в основном следуя [1, 9] – в этих источниках рассматривался случай $E=E'$), получен следующий результат о сходимости метода (3).

**Теорема 1.** *Пусть*

$$\tilde{M}_N^2=\left(\frac{1}{N+1}\sum_{k=0}^{N}M_k^{-2}\right)^{-1},$$

$$R_N^2:=\frac{1}{2}\sum_{e\in E\setminus E'}\left(\tau_e\left(\bar{f}_e^N\right)-\bar{t}_e\right)^2+\frac{1}{2}\sum_{e\in E'}\left(\tilde{t}_e^N-\bar{t}_e\right)^2,$$

$$\left\{\tilde{t}_e^N\right\}_{e\in E'}=\arg\min_{\{t_e\}_{e\in E'}\ge 0}\left\{\underbrace{-\sum_{w\in W}d_wT_w\left(\left\{\tau_e\left(\bar{f}_e^N\right)\right\}_{e\in E\setminus E'},\left\{t_e\right\}_{e\in E'}\right)}_{F\left(\left\{\tau_e\left(\bar{f}_e^N\right)\right\}_{e\in E\setminus E'},\left\{t_e\right\}_{e\in E'}\right)}+\sum_{e\in E'}\bar{f}_e^N\cdot\left(t_e-\bar{t}_e\right)\right\}.$$

*Тогда при*

$$N\ge\frac{2\tilde{M}_N^2R_N^2}{\varepsilon^2}$$

*имеет место неравенство (5) и, как следствие,*

$$0\le\Upsilon\left(\bar{t}^N\right)-\Upsilon_*\le\varepsilon,\quad 0\le\Psi\left(\bar{f}^N\right)-\Psi_*\le\varepsilon. \tag{6}$$

**Доказательство.** Из работ [8, 10], имеем (также в выкладках используется то, что $d\sigma_e^*\left(t_e\right)/dt_e=f_e\Leftrightarrow t_e=\tau_e\left(f_e\right)$, $e\in E\setminus E'$; $\sigma_e^*\left(t_e\right)\ge 0$, $\sigma_e^*\left(t_e^0\right)=\sigma_e^*\left(\bar{t}_e\right)=0$, $e\in E$)

$$\Upsilon\left(\bar{t}^N\right)\le\frac{1}{S_N}\sum_{k=0}^{N}\gamma_k\Upsilon\left(t^k\right)\overset{(*)}{\le}\frac{1}{2S_N}\sum_{k=0}^{N}\gamma_k^2M_k^2+$$

$$+\frac{1}{S_N}\min_{t_e\ge\bar{t}_e,\,e\in E'}\left\{\sum_{k=0}^{N}\gamma_k\left\{F\left(t^k\right)+\left\langle\partial F\left(t^k\right),t-t^k\right\rangle\right\}+S_N\sum_{e\in E}\sigma_e^*\left(t_e\right)+\frac{1}{2}\left\|t-t^0\right\|_2^2\right\}\le\frac{\varepsilon}{2}+$$

$$+\min_t\left\{\frac{1}{S_N}\left\{\sum_{k=0}^{N}\gamma_k\left\{F\left(t^k\right)+\left\langle\partial F\left(t^k\right),t-t^k\right\rangle\right\}+S_N\sum_{e\in E}\sigma_e^*\left(t_e\right)+\frac{1}{2}\left\|t-\bar{t}\right\|_2^2\right\}+\sum_{e\in E'}s_e^N\cdot\left(\bar{t}_e-t_e\right)\right\}\le$$

$$\le\frac{\varepsilon}{2}+\min_{\{t_e\}_{e\in E'}}\min_{\{t_e\}_{e\in E\setminus E'}}\frac{1}{S_N}\left\{\sum_{k=0}^{N}\gamma_k\left\{F\left(t^k\right)+\left\langle\partial F\left(t^k\right),t-t^k\right\rangle\right\}+\right.$$

$$\left.+\left(S_N\sum_{e\in E\setminus E'}\sigma_e^*\left(t_e\right)+\frac{1}{2}\sum_{e\in E\setminus E'}\left(t_e-\bar{t}_e\right)^2\right)+\left(S_N\sum_{e\in E'}\bar{f}_e^N\cdot\left(t_e-\bar{t}_e\right)+\frac{1}{2}\sum_{e\in E'}\left(t_e-\bar{t}_e\right)^2\right)\right\}\le$$

$$\le\frac{\varepsilon}{2}+\min_{\{t_e\}_{e\in E'}\ge 0}\frac{1}{S_N}\left\{\sum_{k=0}^{N}\gamma_k\underbrace{\overbrace{\left\{F\left(t^k\right)+\left\langle\partial F\left(t^k\right),\breve{t}^N-t^k\right\rangle\right\}}^{\breve{t}^N=\left(\left\{\tau_e\left(\bar{f}_e^N\right)\right\}_{e\in E\setminus E'},\left\{t_e\right\}_{e\in E'}\right)}}_{\le F\left(\left\{\tau_e\left(\bar{f}_e^N\right)\right\}_{e\in E\setminus E'},\left\{t_e\right\}_{e\in E'}\right)}+\right.$$

$$\left.+\left(S_N\sum_{e\in E\setminus E'}\sigma_e^*\left(\tau_e\left(\bar{f}_e^N\right)\right)+\frac{1}{2}\sum_{e\in E\setminus E'}\left(\tau_e\left(\bar{f}_e^N\right)-\bar{t}_e\right)^2\right)+\left(S_N\sum_{e\in E'}\bar{f}_e^N\cdot\left(t_e-\bar{t}_e\right)+\frac{1}{2}\sum_{e\in E'}\left(t_e-\bar{t}_e\right)^2\right)\right\}\le$$

$$\le \frac{\varepsilon}{2} + \min_{\{t_e\}_{e\in E'}\ge 0} \left\{ F\left( \left\{ \tau_e\left( \bar{f}_e^N \right) \right\}_{e\in E\setminus E'}, \{t_e\}_{e\in E'} \right) + \right.$$

$$\left. + \sum_{e\in E'} \bar{f}_e^N \cdot \left(t_e - \bar{t}_e\right) + \frac{1}{2S_N} \sum_{e\in E'} \left(t_e - \bar{t}_e\right)^2 \right\} +$$

$$+ \sum_{e\in E\setminus E'} \sigma_e^*\left( \tau_e\left( \bar{f}_e^N \right)\right) + \frac{1}{2S_N} \sum_{e\in E\setminus E'} \left( \tau_e\left( \bar{f}_e^N \right) - \bar{t}_e \right)^2 \le$$

$$\le \frac{\varepsilon}{2} + \min_{\{t_e\}_{e\in E'}\ge 0} \left\{ F\left( \left\{ \tau_e\left( \bar{f}_e^N \right) \right\}_{e\in E\setminus E'}, \{t_e\}_{e\in E'} \right) + \sum_{e\in E\setminus E'} \sigma_e^*\left( \tau_e\left( \bar{f}_e^N \right)\right) + \sum_{e\in E'} \bar{f}_e^N \cdot \left(t_e - \bar{t}_e\right) \right\} +$$

$$+ \frac{1}{2S_N} \sum_{e\in E'} \left(\tilde{t}_e^N - \bar{t}_e\right)^2 + \frac{1}{2S_N} \sum_{e\in E\setminus E'} \left( \tau_e\left( \bar{f}_e^N \right) - \bar{t}_e \right)^2 =$$

$$= \frac{\varepsilon}{2} - \Psi\left(\bar{f}^N\right) + \frac{R_N^2}{S_N} = \frac{\varepsilon}{2} + \frac{\tilde{M}_N^2 R_N^2}{\varepsilon\cdot(N+1)} - \Psi\left(\bar{f}^N\right) \le \varepsilon - \Psi\left(\bar{f}^N\right). \ \square$$

**Следствие 1.** *Пусть* $t^*$ *– решение задачи (2). Положим*

$$R^2 = \frac{1}{2}\left\| t^* - t^0 \right\|_2^2 = \frac{1}{2}\left\| t_* - \bar{t} \right\|_2^2 .$$

*Тогда при*

$$N = \frac{2\tilde{M}_N^2 R^2}{\varepsilon^2}$$

*справедливы неравенства*

$$\frac{1}{2}\left\| t - t^{k+1} \right\|_2^2 \le 2R^2,\ k = 0,\dots,N \tag{7}$$

$$0 \le \Upsilon\left(\bar{t}^N\right) - \Upsilon_* \le \varepsilon . \tag{8}$$

**Доказательство.** Формула (8) – стандартный результат, см., например, [7]. Формула (7) также является достаточно стандартной [10], однако далее приводится схема ее вывода. Из доказательства теоремы 1 имеем для любого $k = 0,\dots,N$

$$0 \le \frac{1}{2}\sum_{l=0}^{k} \gamma_l^2 M_l^2 + \frac{1}{2}\left\| t - t^0 \right\|_2^2 - \frac{1}{2}\left\| t - t^{k+1} \right\|_2^2 .$$

Отсюда следует, что

$$\frac{1}{2}\left\| t - t^{k+1} \right\|_2^2 \le R^2 + \frac{1}{2}\sum_{l=0}^{k} \varepsilon^2 M_l^{-2} \le R^2 + \frac{1}{2}\sum_{k=0}^{N} \varepsilon^2 M_k^{-2} = 2R^2 . \ \square$$

**Замечание 1.** Преимуществом подхода (3), (4), (4’) над подходом п. 3 работы [1] является простота описания (отсутствие необходимости делать рестарты по неизвестным параметрам) и наличие эффективно проверяемого критерия останова (5). К недостаткам стоит отнести вхождение в оценку скорости сходимости плохо контролируемого $R_N^2$, которое может оказаться большим даже в случае, когда $E = E'$. Далее мы опишем другой способ (см. также работы [11, 12], в которых описаны близкие конструкции) восстановления решения прямой задачи (1), отличный от (4), (4’), в части (4’), который в случае $E = E'$ позволяет использовать $R^2$ вместо $R_N^2$.

## 4. Второй способ восстановления решения прямой задачи по решению двойственной

Положим

$$f_e^k \in -\partial_e F\left(t^k\right),\;\; \bar{f}_e^N = \frac{1}{S_N}\sum_{k=0}^{N}\gamma_k f_e^k\,,\; e\in E\,. \tag{9}$$

$$\tilde{R}^2 = \frac{1}{2}\sum_{e\in E'}\left(t_e^* - t_e^0\right)^2 = \frac{1}{2}\sum_{e\in E'}\left(t_e^* - \bar{t}_e\right)^2 .$$

**Теорема 2.** *Пусть*

$$\tilde{R}_N^2 := \frac{1}{2}\sum_{e\in E\setminus E'}\left(\tau_e\left(\bar{f}_e^N\right) - \bar{t}_e\right)^2 + 5\tilde{R}^2 . \tag{10}$$

*Тогда при*

$$N \ge \frac{4\tilde{M}_N^2\tilde{R}_N^2}{\varepsilon^2}$$

*имеют место неравенства*

$$\left|\Upsilon\left(\bar{t}^N\right) - \Upsilon_*\right| \le \varepsilon\,,\; \left|\Psi\left(\bar{f}^N\right) - \Psi_*\right| \le \varepsilon\,.$$

*Более того, также имеют место неравенства*

$$\sqrt{\sum_{e\in E'}\left(\left(\bar{f}_e^N - \bar{f}_e\right)_+\right)^2} \le \tilde{\varepsilon}\,, \tilde{\varepsilon} = \varepsilon/\tilde{R}\,, \tag{11}$$

$$\Psi\left(\bar{f}^N\right) - \Psi_* \le \Upsilon\left(\bar{t}^N\right) + \Psi\left(\bar{f}^N\right) \le \varepsilon\,,$$

*которые можно использовать для критерия останова метода (задавшись парой* $\left(\varepsilon,\tilde{\varepsilon}\right)$*).*

**Доказательство.** Рассуждая аналогично доказательству теоремы 1 (см. также [11, 12]), получим

$$\Upsilon\left(\bar{t}^N\right) \le \frac{1}{S_N}\sum_{k=0}^{N}\gamma_k\Upsilon\left(t^k\right) \overset{(12)}{\le} \frac{1}{2S_N}\sum_{k=0}^{N}\gamma_k^2 M_k^2 +$$

$$+\frac{1}{S_N}\min_{\substack{t_e,\, e\in E\setminus E';\; t_e\ge \bar{t}_e,\, e\in E' \\ \frac{1}{2}\sum_{e\in E'}\left(t_e-\bar{t}_e\right)^2 \le 5\tilde{R}^2}}\left\{\sum_{k=0}^{N}\gamma_k\left\{F\left(t^k\right) + \left\langle \partial F\left(t^k\right), t - t^k\right\rangle\right\} + S_N\sum_{e\in E}\sigma_e^*\left(t_e\right)\right\} + \frac{\tilde{R}_N^2}{S_N} \overset{(13)}{\le}$$

$$\overset{(13)}{\le} \frac{\varepsilon}{2} - \Psi\left(\bar{f}^N\right) - \max_{\substack{t_e\ge \bar{t}_e,\, e\in E' \\ \frac{1}{2}\sum_{e\in E'}\left(t_e-\bar{t}_e\right)^2 \le 5\tilde{R}^2}}\left\{\frac{1}{S_N}\sum_{k=0}^{N}\gamma_k\sum_{e\in E'}\left(f_e^k - \bar{f}_e\right)\left(t_e - \bar{t}_e\right)\right\} + \frac{\tilde{R}_N^2}{S_N} \le$$

$$\le \frac{\varepsilon}{2} - \Psi\left(\bar{f}^N\right) - 3\tilde{R}\sqrt{\sum_{e\in E'}\left(\left(\bar{f}_e^N - \bar{f}_e\right)_+\right)^2} + \frac{\tilde{R}_N^2}{S_N} \le \varepsilon - \Psi\left(\bar{f}^N\right) - 3\tilde{R}\sqrt{\sum_{e\in E'}\left(\left(\bar{f}_e^N - \bar{f}_e\right)_+\right)^2}\,.$$

Неравенство (12) было получено на базе следующего соотношения:

$$\left\{\tau_e\left(\bar{f}_e^N\right)\right\}_{e\in E\setminus E'} = \arg\min_{t_e,\, e\in E\setminus E'}\left\{\frac{1}{S_N}\sum_{k=0}^{N}\gamma_k\left\{F\left(t^k\right) + \left\langle \partial F\left(t^k\right), t - t^k\right\rangle\right\} + \sum_{e\in E}\sigma_e^*\left(t_e\right)\right\}.$$

Неравенство (13) было получено на базе следующих соотношений:

$$\partial F\left(t^k\right) = -f^k,\; F\left(t^k\right) = -\left\langle f^k, t^k\right\rangle,$$

$$\min_{t_e}\left\{-f_e^k t_e + \sigma_e^*\left(t_e\right)\right\} = -\sigma_e\left(f_e^k\right),\; e\in E\setminus E'\,,$$

$$-f_e^k t_e + \sigma_e^*(t_e) = -\left(f_e^k - \bar{f}_e\right)\left(t_e - \bar{t}_e\right) - f_e^k \bar{t}_e = -\left(f_e^k - \bar{f}_e\right)\left(t_e - \bar{t}_e\right) - \sigma_e\left(f_e^k\right),\ e \in E',$$

$$-\frac{1}{S_N}\sum_{k=0}^{N}\gamma_k \sum_{e\in E}\sigma_e\left(f_e^k\right) \le -\Psi\left(\bar{f}^N\right).$$

Таким образом,

$$\Upsilon\left(\bar{t}^N\right) + \Psi\left(\bar{f}^N\right) + 3\tilde{R}\sqrt{\sum_{e\in E'}\left(\left(\bar{f}_e^N - \bar{f}_e\right)_+\right)^2} \le \varepsilon .$$

Повторяя рассуждения п. 6.11 [13] и п. 3 [12], получим искомые неравенства. □

## 5. Заключительные замечания

В данном пункте в формате замечаний обсуждаются полученные в статье результаты и кратко отмечаются результаты численных экспериментов.

**Замечание 2.** О преимуществе использования формулы (9) вместо (4), (4’) в случае $E = E'$ написано в замечании 1. Сейчас отметим возникающие при таком подходе недостатки: 1) возможность нарушения ограничения $f_e \le \bar{f}_e$, $e \in E'$ в прямой задаче, 2) отсутствие левых неравенств в двойных неравенствах (5), (6).

**Замечание 3.** Формулы (4), (9) вынужденно (в случае (4)) или осознано (в случае (9) при $e \in E'$) восстанавливают решение прямой задачи исходя из “модели” – явной формулы, связывающей прямые и двойственные переменные. Наличие таких переменных неизбежно приводит к возникновению в оценках зазора двойственности трудно контролируемых размеров решений вспомогательных задач. При этом в случае, когда наличие модели сопряжено с каким-то ограничением в прямой задаче (формулы (9) при $e \in E'$, и ограничение в прямой задаче $f_e \le \bar{f}_e$, $e \in E'$), допускается нарушение этих ограничений, которые необходимо контролировать (11). Зато по этим переменным имеется полный контроль соответствующих этим переменным частей оценок зазора двойственности (10). Подход (4’) связанный с наличием ограничений в решаемой задаче ($t_e \ge \bar{t}_e$, $e \in E'$) также приводит к возникновению в оценках зазора двойственности трудно контролируемых размеров решений вспомогательных задач, однако уже не приводит к нарушению никаких ограничений в самой задаче и сопряженной к ней (в нашем случае исходная задача – двойственная (2), а сопряженная к ней – прямая (1)). Эти два прямо-двойственных подхода дополняются другим прямо-двойственным подходом, в котором шаги осуществляются “по функционалу”, если не нарушены или слабо нарушены ограничения в задаче, и “по нарушенному ограничению” в противном случае. Подробнее об этом см., например, в работах [13, 14].

**Замечание 4.** Оба описанных подхода можно распространить, сохраняя вид формул восстановления и структуру рассуждений, на практически произвольные пары прямая / двойственная задача, поскольку выбранный нами пример пары взаимно-сопряженных задач (1), (2) и так содержал в себе практически все основные нюансы, которые могут возникать при таких рассуждениях. Более того, вместо МЗС можно было бы использовать любой другой прямо-двойственный метод. Например, композитный универсальный градиентный метод из [15].

**Замечание 5.** Возможность эффективного распараллеливания предложенных методов связана с возможностью эффективного вычисления самой затратной части шага описанного итерационного метода: расчет элемента субдифференциала $\partial F\left(t^k\right)$ (см. формулы (3), (4), (9) в которых этот субдифференциал используется):

$$\partial F(t) = -\sum_{i\in O}\sum_{j\in D:(i,j)\in W} d_{ij}\partial T_{ij}(t).$$

Вычисление $\left\{\partial T_{ij}(t)\right\}_{j\in D:(i,j)\in W}$ может быть осуществлено алгоритмом Дейкстры [16] (и его более современными аналогами [17]) за $O(n\ln n)$. При этом под $\partial T_{ij}(t)$ можно понимать описание произвольного (если их несколько) кратчайшего путь из вершины $i$ в вершину $j$ на графе $\Gamma$, ребра которого взвешены вектором $t=\{t_e\}_{e\in E}$. Под “описанием” понимается $\left[\partial T_{ij}(t)\right]_e=1$, если $e$ попало в кратчайший путь и $\left[\partial T_{ij}(t)\right]_e=0$ иначе. Таким образом, вычисление $\partial F(t)$ может быть распараллелено на $S$ процессорах.

**Замечание 6.** Строго говоря, нужно найти вектор $\partial F(t)$, а не кратчайшие пути. Чтобы получить вектор $\partial F(t)$ за $O(Sn\ln n)$ стоит для каждого из $S$ источников построить (например, алгоритмом Дейкстры) соответствующее дерево кратчайших путей. Исходя из принципа динамического программирования “часть кратчайшего пути сама будет кратчайшим путем” несложно понять, что получится именно дерево, с корнем в рассматриваемом источнике. Это можно сделать для одного источника за $O(n\ln n)$ [16, 17]. Однако, главное, правильно взвешивать ребра (их не больше $O(n)$) такого дерева, чтобы за один проход этого дерева можно было восстановить вклад (по всем ребрам) соответствующего источника в общий вектор $\partial F(t)$. Ребро должно иметь вес равный сумме всех проходящих через него корреспонденций с заданным источником (корнем дерева). Имея значения соответствующих корреспонденций (их также не больше $O(n)$) за один обратный проход (то есть с листьев к корню) такого дерева можно осуществить необходимое взвешивание (с затратами не более $O(n)$). Делается это по правилу: вес ребра равен сумме корреспонденции (возможно, равной нулю), в соответствующую вершину, в которую ребро входит и сумме весов всех ребер (если таковые имеются), выходящих из упомянутой вершины.

Студентами 6-го курса ФУПМ и ФАКИ МФТИ были проведены разнообразные численные эксперименты [18, 19] с описанными в данной работе и работе [1] алгоритмам. Данные для экспериментов брались с открытого ресурса [20]. Ниже собраны основные выводы, сделанные из этих экспериментов:

- если в смешанной модели все ребра Бэкмановского типа, то лучше использовать метод условного градиента [1, 2] (Франк–Вульфа);
- способы, описанные в данной статье и в работе [1], базирующиеся на решении двойственной задачи прямо-двойственным методом зеркального спуска работают согласно полученным оценкам (то есть полученные оценки подтверждаются практикой, т.е. быстрее, чем предписывают полученные оценки, разработанные методы не сходятся);
- применять двойственные подходы, в основе которых лежит зеркальный спуск (метод двойственных усреднений [10]) к реальным городам оправдано, только если ограничиваться низкой относительной точностью (не точнее 5 %) – для мегаполисов (десятки тысяч ребер) все эти методы будут работать часами (не параллельный вариант) для достижения более высокой точности;
- использование языков высокого уровня (типа Python), приводит к потери эффективности для мегаполисов (вплоть до потери одного порядка) по сравнению с языками более низкого уровня (типа C++).

Авторы выражают благодарность А.С. Аникину, П.Е. Двуреченскому, М.Б. Кубентаевой, А.А. Шананину, С.В. Шпирко и студентам МФТИ за помощь в работе.

Работа А.В. Гасникова и Е.В. Гасниковой в пп. 2, 3 выполнена при поддержке гранта Президента РФ МД 1320.2018.1. Работа А.В. Гасникова и Ю.Е. Нестерова в пп. 4, 5 выполнена при поддержке гранта РНФ 17-11-01027.

## Литература